\documentclass[amscd,syntonly,amssymb]{amsart}

\usepackage{epsfig}
\theoremstyle{plain}
\newtheorem{Thm}{Theorem}
\newtheorem{Prop}[Thm]{Proposition}
\newtheorem{Cor}[Thm]{Corollary}
\newtheorem{Lem}[Thm]{Lemma}

 \theoremstyle{definition}
\theoremstyle{remark}

\numberwithin{equation}{section}

\begin{document}
 \title{A characteristic number of bundles determined by mass linear pairs}

 \author{ ANDR\'{E}S   VI\~{N}A}
\address{Departamento de F\'{i}sica. Universidad de Oviedo.   Avda Calvo
 Sotelo.     33007 Oviedo. Spain. }
 \email{vina@uniovi.es}
\thanks{This work has been partially supported by Ministerio de Ciencia y
Tecnolog\'{\i}a, grant MAT2007-65097-C02-02}

\begin{abstract}

 Let $\Delta$ be a Delzant polytope in ${\mathbb R}^n$  and ${\bf b}\in{\mathbb
 Z}^n$.
  Let $E$ denote the symplectic fibration over $S^2$ determined by the pair $(\Delta,\,{\bf b})$.
We prove the equivalence between the
  fact that  $(\Delta,\,{\bf b})$ is a mass linear pair (D. McDuff, S. Tolman,
   {\em Polytopes with mass linear functions, part I.} {\tt
   arXiv:0807.0900 [math.SG]})
 and the vanishing
 of a characteristic number of  $E$ in the following cases: When
 $\Delta$ is  a $\Delta_{n-1}$ bundle over
 $\Delta_1$;
  when $\Delta$ is
the polytope associated with the one point blow up of ${\mathbb
C}P^n$;  and  when $\Delta$ is the  polytope associated with a
Hirzebruch
 surface.

\end{abstract}

 \maketitle

   \smallskip


 \medskip

\centerline{}

\bigskip
\section {Introduction}

 Let $T$ be the torus $(U(1))^n$ and $\Delta=\Delta({\bf n},k)$ the
  polytope
in ${\mathfrak t}^*$ with $m$ facets defined by
\begin{equation}\label{DeltP}
\Delta({\bf n},\,k)=\bigcap_{j=1}^m\,\{ x\in {\mathfrak
t}^*\,:\,\langle x,{\bf n}_j\rangle\leq k_j \},
\end{equation}
where $k_j\in{\mathbb R}$ and the ${\bf n}_j\in {\mathfrak t}$ are
the outward conormals to the facets. The facet defined by the
equation $\langle x,\,{\bf n}_j\rangle=k_j$ will be denoted $F_j$,
and we put $\text{Cm}(\Delta)$ for the
 mass center of the polytope $\Delta$.

In \cite{M-T3} is defined the chamber  ${\mathcal C}_{\Delta}$ of
$\Delta$ as the set of $k'\in{\Bbb R}^m$ such that the polytope
$\Delta':=\Delta({\bf n},\,k')$ is analogous to $\Delta$; that is,
 the intersection $\cap_{j\in J} F_j$ is nonempty
iff $\cap_{j\in J} F'_j\ne \emptyset$ for any
$J\subset\{1,\dots,m\}$.
When we consider only polytopes which belong to the chamber of a
fixed polytope we delete the ${\bf n}$ in the notation introduced
in (\ref{DeltP}).

 Further McDuff and
Tolman introduced the concept of mass linear pair: Given the
polytope
  $\Delta$ and ${\bf b}\in{\mathfrak t}$, the pair
  $( \Delta,\, {\bf b})$ is mass linear if the map
$$k\in{\Bbb R}^m\mapsto \langle \text{Cm}(\Delta(k)),\,{\bf b}\rangle\in{\Bbb R}$$
is linear on
 ${\mathcal C}_{\Delta}$.

Let $(N,\Omega)$ be a closed connected symplectic $2n$-manifold.
By $\text{Ham}(N,\Omega)$ we denote the Hamiltonian group of
$(N,\Omega)$ \cite{Mc-S}, \cite{lP01}. If $\psi$ is a loop in
$\text{Ham}(N,\Omega)$, then $\psi$ determines a Hamiltonian fibre
bundle $E\to S^2$ with standard fibre $N$ via the clutching
construction. In \cite{L-M-P} various characteristic numbers for
the fibre bundle $E$   are defined. These numbers give rise to
topological invariants of the loop $\psi$. In this note we will
consider only the following characteristic number
 \begin{equation}\label{Ipsi}
 I({\psi}):=\int_Ec_1(VTE)\,c^n,
 \end{equation}
 where $VTE$ is the vertical tangent bundle of  $E$ and $c\in H^2(E,\,{\Bbb R})$ is the
 coupling class of the fibration $E\to S^2$ \cite{G-L-S}, \cite{Mc-S}.
 $I({\psi})$ depends only on the homotopy class of the loop $\psi$.
 Moreover the map
 \begin{equation}\label{Ihomo}
  I:\psi\in\pi_1(\text{Ham}(N))\mapsto I({\psi})\in{\Bbb R}
   \end{equation}
 is an ${\Bbb R}$-valued group homomorphism \cite{L-M-P}.

Let us suppose that $\Delta$ is a Delzant polytope. We shall
denote by $(M_{\Delta},\,\omega_{\Delta},\,\mu_{\Delta})$ the
toric manifold determined by $\Delta$
 ($\mu_{\Delta}:M\to{\mathfrak t}^*$  being the corresponding moment map).
 Given ${\bf b}$, an element in the integer lattice of ${\mathfrak
t}$, we shall write $\psi_{\bf b}$ for the loop of Hamiltonian
diffeomorphisms of $(M_{\Delta},\,\omega_{\Delta})$ defined by
${\bf b}$.
 The   bundle with fibre $M_{\Delta}$ determined by
$\psi_{\bf b}$ will be denoted  $E_{\Delta,\,{\bf b}}$, and we
will let $I(\Delta;\,{\bf b})$ for the characteristic number
$I(\psi_{\bf b})$.
 When we consider only
polytopes
 in the chamber of a given polytope, we will write $I(k;\,{\bf b})$
 instead of $I(\Delta(k);\,{\bf b})$ for $k$ in this chamber.

In Section 2 we   study the characteristic number $I(k;\,{\bf
b})$, when $(\Delta,\,{\bf b})$ is linear
 pair and $k$ varies in the chamber of $\Delta$, and
 we    prove that $I(k;\,{\bf b})$ is a homogeneous polynomial of ten $k_j$
 (see Proposition \ref{polyno}).

 In Section 3 we consider the case when    the polytope $\Delta$ is a $\Delta_p$ bundle over $\Delta_1$
with $p>1$ \cite{M-T3}.  Then   $M_{\Delta}$ is a
$2(p+1)$-dimensional manifold diffeomorphic to the total space of
the fibre bundle ${\mathbb P}(L_1\oplus\dots\oplus
L_p\oplus{\mathbb C})\to{\mathbb C}P^1 $, where each $L_j$ is a
holomorphic line bundle over ${\mathbb C}P^1$. Given ${\bf
b}\in{\mathbb Z}^{p+1}$, we prove that $I(k;\,{\bf b})$ is the
product of two factors, ${\mathcal K}$ and ${\mathcal Z}({\bf
b})$, such that the first one is independent of ${\bf b}$  and the
second one is independent of $k\in{\mathcal C}_\Delta$ (Theorem
\ref{ThmHatb+Dotb}). We also prove the equivalence between the
vanishing of ${\mathcal Z}({\bf b})$ and the fact that
$(\Delta,\,{\bf b})$ is a mass linear pair (Theorem
\ref{ThmMasLinear}). As a consequence we deduce that a necessary
and sufficient condition for the vanishing of $I(k;\,{\bf b})$ on
${\mathcal C}_{\Delta}$ is that $(\Delta,\,{\bf b})$ be a mass
linear pair (see Theorem \ref{vanish}).
 From this theorem
 we will deduce that $\psi_{\bf b}$ generates an infinite cyclic
  subgroup in $\pi_1({\rm Ham}(M_{\Delta}))$, if the pair $(\Delta,\,{\bf
  b})$ is not mass linear (Proposition \ref{Proppi}).

The case when the polytope $\Delta$ is the one associated with a
Hirzebruch surface is considered in Subsection \ref{SubsectHirz}.
Using calculations carried out in \cite{V},  we will prove the
equivalence between the vanishing of $I(\Delta;\,{\bf b})$ and the
fact that $(\Delta,\,{\bf b})$ is a mass linear pair (Theorem
\ref{ThmHirz}). In a Remark we give a general proof, based in
general properties of the characteristic number $I$, of the
implication: $(\Delta,{\bf b})$ is a mass linear pair
$\Longrightarrow\,$ $I(k;\,{\bf b})$ vanishes on the chamber of
$\Delta$. The arguments developed in this proof are applicable to
 other polytopes; for example to the polytope associated
with the one point blow up of ${\mathbb C}P^n$ (see Proposition
\ref{Propfiniteorder}).

 In Subsection \ref{SubsectCP2} we consider the polytope $\Delta$
 associated to  the manifold one point blow up of ${\mathbb C}P^n$.
 We also prove the equivalence: $I(k;\,{\bf
 b})=0$ for all $k\in{\mathcal C}_{\Delta}$ $\,\Longleftrightarrow\, (\Delta,\,{\bf b})$ is a mass linear pair (see
 Theorem \ref{Thmequiv}).
 As a consequence we deduce a simple sufficient condition for  $\psi_{\bf b}$ to generate an
infinite cyclic subgroup in
 $\pi_1({\rm Ham}(M_{\Delta}))$ (Proposition \ref{Propblowup}).

In summary, we prove
 the equivalence between the vanishing of
 $I(k;\,{\bf b})$ for all $k\in{\mathcal C}_{\Delta}$ and the
 property of  $(\Delta,\,{\bf b})$ being a mass linear pair,
 in the
following cases: When $\Delta$ is  a $\Delta_p$ bundle over
$\Delta_1$, when $\Delta$ is the trapezoid associated to a
Hirzebruch  surface, when $\Delta$ is the truncated simplex
associated to the one point blow up of ${\mathbb C}P^n$.

In the proof of the mentioned  results plays a crucial role a
formula for
 the characteristic number $I({\psi_{\bf b}})$ obtained in
\cite{V2}. This formula gives $I({\psi_{\bf b}})$
  in terms of the integrals, on the facets of
the polytope, of the normalized Hamiltonian corresponding to the
loop $\psi_{\bf b}$ (see (\ref{Ibr})).

Let   $(\Delta,\,{\bf b})$ be a pair consisting of a Delzant
polytope in ${\mathfrak t}^*$
  and an element in the integer lattice of ${\mathfrak t}$.
In view the above results one is tempted to conjecture the
equivalence between the following statements

a)  $I(k;\,{\bf b})=0$ for all $k\in{\mathcal C}_{\Delta}$.

b)   $(\Delta,\,{\bf b})$ is a  mass linear pair.

 We think that a possible proof of this conjecture using formula (\ref{Ibr})
will probably involve general properties, valid for all Delzant
polytopes, about $\text{Cm}(\Delta)$ and the mass center of the
facets $F_j$.




\medskip

{\it Acknowledgements.} I thank Dusa McDuff and Susan Tolman for
sending me a working draft of the paper \cite{M-T3} and for
comments.

\smallskip


\section {A characteristic number}\label{SecChar}

Let us suppose that the polytope $\Delta$  defined in
(\ref{DeltP}) is a Delzant polytope
 in ${\mathfrak t}^*$. If ${\bf b}$ is in
the integer lattice of ${\mathfrak t}$,  an expression for the
value of $I({\psi_{\bf b}})$ in terms of integrals of the
Hamiltonian function has been obtained in Section 4 of \cite{V2}
\begin{equation}\label{Ibr}
I(\Delta;\,{\bf b}):=I({\psi_{\bf
b}})=-n\sum_{j=1}^m\int_{D_j}f\,(\omega_{\Delta})^{n-1},
\end{equation}
  where $D_j:=\mu_{\Delta}^{-1}(F_j)$ is oriented by the restriction of $\omega_{\Delta}$,  and $f$ being the normalized Hamiltonian of the corresponding circle
action; that is,
 $$f=\langle \mu_{\Delta},\,{\bf b}\rangle+ \text{constant}\;\;\;\text{and}\;\;\;
\int_{M_{\Delta}}f\,(\omega_{\Delta})^n=0.$$

That is,
\begin{equation}\label{Ib}
I(\Delta;\,{\bf b})=n\sum_{j=1}^m\Big( \langle
\text{Cm}(\Delta),\,{\bf
b}\rangle\int_{D_j}(\omega_{\Delta})^{n-1}-\int_{D_j}\langle\mu_{\Delta},\,{\bf
b}\rangle \,(\omega_{\Delta})^{n-1} \Big),
\end{equation}
where
\begin{equation}\label{Cm}
\langle \text{Cm}(\Delta),\,{\bf b}\rangle= \frac{\int_M
\langle\mu_{\Delta},\,{\bf b}\rangle\,(\omega_{\Delta})^n}{\int_M
(\omega_{\Delta})^n}.
\end{equation}

If $\Delta=\Delta({\bf n},\,k)$ we consider the polytope
 $\Delta'=\Delta({\bf n},\,k')$ obtained from
$\Delta$ by the translation defined by a vector $a$ of ${\mathfrak
t}^*$. As we said, we write $I(k;\,{\bf b})$ and $I(k';\,{\bf b})$
for the corresponding characteristic numbers. According to the
construction of the respective toric manifolds (see \cite{Gui}),
$$M_{\Delta'}=M_{\Delta},\;\;\;
\omega_{\Delta'}=\omega_{\Delta},\;\;\;
\mu_{\Delta'}=\mu_{\Delta}+a.$$
  But the {\em normalized}
Hamiltonians $f$ and $f'$ corresponding to the action of ${\bf b}$
on $M_{\Delta}$ and $M_{\Delta'}$ are equal. Thus it follows from
(\ref{Ibr}) that $I(k;\,{\bf b})=I(k';\,{\bf b})$. More precisely,
we have the evident proposition

\begin{Prop}\label{Propinv}
If $a$ is an arbitrary vector of ${\mathfrak t}^*$, then
$I(k;\,{\bf b})=I(k';\,{\bf b})$, for $k'_i=k_i+\langle a,{\bf
n}_i\rangle$, $\;i=1,\dots, m$.

\end{Prop}

\smallskip

Following \cite{Gui} we recall some points of the construction of
$(M_{\Delta},\,\omega_{\Delta},\,\mu_{\Delta})$ from the polytope
$\Delta$ defined by (\ref{DeltP}), in order to study the value of
the integrals that appear in (\ref{Ib}) and (\ref{Cm}). Given
$\Delta$, we put $r:=m-n$ and $\tilde T:=(S^1)^r$. The ${\bf n}_i$
determine weights $w_j\in\tilde{\mathfrak t}^*$, $j=1,\dots,m$ for
a $\tilde T$-action on ${\Bbb C}^m$. Then
moment map for this action is
$$J:z\in{\Bbb C}^m\mapsto J(z)=\pi\sum_{j=1}^m|z_j|^2w_j\in \tilde{\mathfrak
t}^*.$$
 The $k_i$ define a regular value $\sigma$ for $J$, and the
 manifold $M$ is the following orbit space
 \begin{equation}\label{DefiM}
 M_{\Delta}=\{z\in {\Bbb C}^m\,:\,\pi\sum_{j=1}
^m|z_{j}|^2w_{j}=\sigma\}/{\tilde T},
 \end{equation}
 where the relation defined by ${\tilde T}$ is
 \begin{equation}\label{DefiAc}(z_j)\simeq(z'_j)\;\;\text{iff
there is}\;\,{\bf y}\in \tilde {\mathfrak t}\;\, \text{such
that}\;\, z'_j=z_j\,e^{2\pi i\langle w_j,{\bf y} \rangle}\;
\,\text{for}\;\, j=1,\dots, m.
 \end{equation}

Identifying $\tilde{\mathfrak t}^*$ with ${\Bbb R}^r$,
$\sigma=(\sigma_1,\dots,\sigma_r)$ and each $\sigma_a$ is a linear
combination of the $k_j$'s.

 After a possible change in numeration
 of the facets, we can assume that $F_1,\dots,F_n$ intersect at a
vertex of $\Delta$. If we write $z_j=\rho_je^{i\theta_j}$, then
the symplectic form can be written on $\{[z]\in M\,:\, z_i\ne 0,\,
\forall i \}$
 \begin{equation}\label{omeg}
  \omega_{\Delta}=(1/2)\sum_{i=1}^n d\rho_i^2\wedge
d\varphi_i,
 \end{equation}
  with $\varphi_i$ an angular variables, linear combination of the $\theta_j$'s.

The action of $T=(S^1)^n$ on $M$
$$(\alpha_1,\dots,\alpha_n)[z_1,\dots,
z_m]:=[\alpha_1z_1,\dots,\alpha_nz_n,z_{n+1},\dots,z_m]$$
 gives  $M$ a structure of  toric manifold.
Identifying ${\mathfrak t}^*$ with ${\mathbb R}^{n}$,
 the moment map
 $\mu_{\Delta}:M_{\Delta}\to{\mathfrak t}={\Bbb R}^n$ is defined by
 $$\mu_{\Delta}([z])=\pi(\rho^2_1,\dots,\rho_n^2)+(d_1,\dots,d_n),$$
where the constants $d_i$  are linear combinations of the $k_j$'s
and
 \begin{equation}\label{immu}
  \text{im}\,\mu_{\Delta}=\Delta.
 \end{equation}

 By Proposition \ref{Propinv} we can assume that all $d_j$ are
 zero in the determination of $I(k;\,{\bf b})$.

 We write $x_i:=\pi\rho_i^2$, then
 $$\int_{M_{\Delta}}(\omega_{\Delta})^n=
n! \int_{\Delta} dx_1\dots dx_n,\;\; \;\;
  \int_{M_{\Delta}}\langle\mu_{\Delta},\,{\bf b}\rangle(\omega_{\Delta})^n=
 n!\int_{\Delta}\sum_{i=1}^n b_ix_i\, dx_1\dots dx_n.$$

The following Lemma is useful to evaluate some integrals which
will appear henceforth.
\begin{Lem}\label{Lemmaint}
If
 $$S_n(\tau):=\Big\{(x_1,\dots,x_n)\in{\mathbb
R}^n\,\Big|\,\sum_{i=1}^nx_i\leq\tau,\;\;\;0\leq x_j, \;\forall
j\Big\},$$ then
$$\int_{S_n(\tau)}f(x_1,\dots,x_n)\,dx_1\dots dx_n=\begin{cases}

 \frac{\tau^{n}}{n!}\, ,\; \text{if}\;\;\; f=1 \\

                                         \\

c\frac{\tau^{n+c}}{(n+c)!}\, ,\; \text{if}\;\;\; f=x_i^c,\; c=1,2 \\
              \\
\frac{\tau^{n+2}}{(n+2)!}\, ,\; \text{if}\;\;\; f=x_i x_j,\;  i\ne
j.
\end{cases}$$

\end{Lem}

More general, if $c_1,\dots,c_n\in{\mathbb R}_{>0}$, we put
$$S_n(c,\tau):=\Big\{(x_1,\dots,x_n)\in{\mathbb
R}^n\,\Big|\,\sum_{i=1}^n c_ix_i\leq\tau,\;\;\;0\leq x_j,
\;\forall j\Big\},$$
 then
 \begin{equation}\label{Intsnctau}
  \int_{S_n(c,\tau)}dx_1\dots
 dx_n=\frac{1}{n!}\prod_{i=1}^n\frac{\tau}{c_i},\;\;\; \; \int_{S_n(c,\tau)}x_j\,dx_1\dots
 dx_n=\frac{1}{(n+1)!}\frac{\tau}{c_j}\prod_{i=1}^n\frac{\tau}{c_i}
  \end{equation}
Thus, in the particular case  that $\Delta=S_n(c,\tau)$, then
$\int_{M_{\Delta}}(\omega_{\Delta})^n$ is a monomial of degree $n$
in $\tau$, and $\int_{M_{\Delta}}\langle \mu_{\Delta},\,{\bf
b}\rangle(\omega_{\Delta})^n$ is a monomial of degree $n+1$.

\smallskip

We return to the general case in which  $\Delta$ is the polytope
defined in (\ref{DeltP}). Its
 vertices   are the solutions to
$$\langle x,\,{\bf n}_{j_a}\rangle=k_{j_a},\;\, a=1,\dots,n,$$
hence the coordinates of any vertex of $\Delta$ are linear
combinations of the $k_j$.

A hyperplane in ${\mathbb R}^n$ through a vertex
$(x^0_1,\dots,x_n^0)$ of $\Delta$ is given by an equation of the
form
 \begin{equation}\label{langlex}
  \langle x,\,{\bf n}\rangle =\langle x^0,\,{\bf
n}\rangle=:\kappa.
 \end{equation}
 So $\kappa$ is a linear combination of the $k_j$.

 By drawing  hyperplanes   through vertices of $\Delta$ we can obtain
 a family $\{\,_{\beta}S\}$ of subsets of $\Delta$ such that:

 a) Each $_{\beta}S$ is the transformed of a simplex
 $S_n(b,\tau)$ by an element of the
 group of Euclidean motions in ${\mathbb R}^n$.

 b) For $\alpha\ne\beta$, $\, \,_{\alpha}S\cap\,_{\beta}S$  is a
 subset of the border of $_{\alpha}S$.

 c) $\bigcup_{\beta}\,_{\beta}S=\Delta$.

So, by construction, each facet of $_{\beta}S$ is of the form
(\ref{langlex}) with $\kappa$ linear combination of the $k_j$.

On the other hand the hyperplane $\pi$, $\langle x,\,{\bf
n}\rangle=\kappa$, is transformed by an element of $SO(n)$ in an
hyperplane $\langle x,\,{\bf n'}\rangle=\kappa$.  If ${\mathcal
T}$ is a translation in ${\mathbb R}^n$ which applies
$S_n(b,\tau)$ onto $_{\beta}S$, then this transformation maps
$(0,\dots, 0)$ in a vertex $a=(a_1,\dots,a_n)$ of $_{\beta}S$. So
the translation ${\mathcal T}$ transforms $\pi$ in $\langle
x,\,{\bf n}\rangle=\kappa+\langle a,\,{\bf n}\rangle=:\kappa'$. As
the $a_j$ are linear combinations of the $k_j$, so is $\kappa'$.
Hence any element of the group of Euclidean motions in ${\mathbb
R}^n$ which maps  $S_n(b,\tau)$ onto $_{\beta}S$ transforms the
hyperplane $\pi$ through a vertex of $\Delta$ in an hyperplane
$\langle x,\,{\bf n'}\rangle=\kappa'$ with $\kappa'$ a linear
combination of the $k_j$. In particular, $\tau$ is a linear
combination of the $k_j$, and   by (\ref{Intsnctau})
$$\int_{_\beta S}dx_1\dots dx_n=\int_{S_{n}(b,\tau)}dx_1\dots
dx_n$$
 is a monomial of degree $n$ of a linear combination of the $k_j$.
 Thus,
 $$\int_{M} (\omega_{\Delta})^n=\sum_{\beta}\int_{_\beta
 S}dx_1\dots dx_n,$$
is a homogeneous polynomial of degree $n$ of the $k_j$.

Similarly
$$\int_{M_{\Delta}}\langle\mu_{\Delta},\,{\bf
b}\rangle(\omega_{\Delta})^n$$
 is a homogeneous polynomial of degree $n+1$ of the $k_j$.
 Analogous results hold for $\int_{D_j}(\omega_{\Delta})^{n-1}$ and
$\int_{D_j}\langle\mu_{\Delta},\,{\bf
b}\rangle(\omega_{\Delta})^{n-1}$.

 It follows from (\ref{Ib}), (\ref{Cm}) together with the
 preceding argument
   the following
 proposition
 \begin{Prop}\label{HomPr}
 Given a Delzant polytope $ \Delta$, if ${\bf b}$ belongs to
 the integer lattice of ${\mathfrak t}$, then
 $I(k;\,{\bf b})$
is a rational function of the $k_i$ for $k\in{\mathcal
C}_{\Delta}$.
\end{Prop}

 Analogously we have

\begin{Prop}\label{polyno}
If $(\Delta,\,{\bf b})$ is mass linear pair, then $I(k;\,{\bf b})$
is a homogeneous polynomial in the $k_i$ of degree $n$, when
$k\in{\mathcal C}_{\Delta}$.
\end{Prop}

\begin{Prop}\label{PropDivisible} If for the polytope $\Delta$ defined in (\ref{DeltP}) ${\bf
n}_a=-{\bf n}_i$, with $i\ne a$, and $(\Delta,\,{\bf b})$ is a
mass linear pair, then  $I(k;\,{\bf b})$ is a polynomial divisible
by $k_a-k_i$.
\end{Prop}

{\it Proof.}  Given $k\in{\mathcal C}_{\Delta}$, maintaining $k_a$
fixed we vary $k$ so that $k\in{\mathcal C}_{\Delta}$ and
$k_i\to k_a$. In the limit $\Delta$ collapses in the facet $F_a$.
As $(\Delta,\,{\bf b})$ is mass linear
 \begin{equation}\label{limitF_a}
  \lim_{k_i\to k_a}\langle\text{Cm}(\Delta),\,{\bf b}\rangle =
 \langle\text{Cm}(F_a),\,{\bf b}\rangle,
 \end{equation}
  where
  \begin{equation}\label{CmF_a}
  \langle\text{Cm}(F_a),\,{\bf b}\rangle=
   \frac{\int_{D_a}\langle\mu_{\Delta},\,{\bf
   b}\rangle\,(\omega_{\Delta})^{n-1}}{
   \int_{D_a} (\omega_{\Delta})^{n-1}}.
  \end{equation}

 We write (\ref{Ib}) as
 \begin{equation}\label{IbmathE}
  I(k;\,{\bf b})=n\sum_{j=1}^{m}{\mathcal E}_j,
  \end{equation}
   with
$${\mathcal E}_j=\langle\text{Cm}(\Delta),\,{\bf
b}\rangle\int_{D_j}(\omega_{\Delta})^{n-1}-
\int_{D_j}\langle\mu_{\Delta},\,{\bf
b}\rangle\,(\omega_{\Delta})^{n-1}.$$

In this limit process the facet $F_a$ remains unchanged, so by
(\ref{limitF_a}) and (\ref{CmF_a})
 \begin{equation}\label{LimtE_a}
 \lim_{k_i\to k_a}{\mathcal E}_a=
  \langle\text{Cm}(F_a),\,{\bf
b}\rangle\int_{D_a}(\omega_{\Delta})^{n-1}-
\int_{D_a}\langle\mu_{\Delta},\,{\bf
b}\rangle\,(\omega_{\Delta})^{n-1}=0.
 \end{equation}

On the other hand, the facets $F_j$, with $j\ne a$ give rise, in
the limit $k_i\to k_a$, to a subdivision of $F_a$
(in Remark after Theorem \ref{ThmHirz} is detailed  this
subdivision in a
  particular case).
  Hence
\begin{align}\label{Sum_jnea}
\lim_{k_i\to k_a}\sum_{j\ne
a}\int_{D_j}(\omega_{\Delta})^{n-1}=&\int_{D_a}(\omega_{\Delta})^{n-1},
\\ \label{alignsecond}
 \lim_{k_i\to k_a}\sum_{j\ne
a}\int_{D_j}\langle\mu_{\Delta},\,{\bf b}\rangle\,
 (\omega_{\Delta})^{n-1}=&\int_{D_a}\langle\mu_{\Delta},\,{\bf
b}\rangle\,(\omega_{\Delta})^{n-1}.
\end{align}

From (\ref{limitF_a}), (\ref{Sum_jnea}), (\ref{alignsecond})  and
(\ref{CmF_a}) it follows
\begin{equation}\label{limSumE}
\lim_{k_i\to k_a} \sum_{j\ne a}{\mathcal E}_j=
 \langle\text{Cm}(F_a),\,{\bf
b}\rangle\int_{D_a}(\omega_{\Delta})^{n-1}-
\int_{D_a}\langle\mu_{\Delta},\,{\bf
b}\rangle(\omega_{\Delta})^{n-1}=0.
\end{equation}

We deduce from (\ref{LimtE_a}) and (\ref{limSumE}) that
$\lim_{k_i\to k_a} I(k;\,{\bf
 b})=0$. That is, $I(k;\,{\bf b})$ is divisible by $k_i-k_a$.
 \qed

\smallskip


\section{$\Delta_p$ bundle over $\Delta_1$}

 Given the integer $p>1$, following \cite{M-T3} we consider the following vectors in ${\mathbb R}^{p+1}$
  \begin{equation}\label{definition-ni}
   {\bf n}_i=-e_i,\; i=1,\dots,p,\;\;{\bf n}_{p+1}=\sum_{i=1}^p
  e_i,\;\; {\bf n}_{p+2}=-e_{p+1},\;\; {\bf
 n}_{p+3}=e_{p+1}-\sum_{i=1}^pa_ie_i,
  \end{equation}
  where $e_1,\dots, e_{p+1}$
is the standard basis of ${\mathbb R}^{p+1}$ and $a_i\in{\mathbb
Z}.$ We write
$${\bf a}:=(a_1,\dots,a_p)\in{\mathbb Z}^{p},\;\;A:=\sum_{i=1}^p
a_i.$$

 Let  $\lambda,\,\tau$  be real positive numbers with $\lambda+a_i>0,$
for $i=1,\dots,p$. We will consider the polytope $\Delta$ in
$({\mathbb R}^{p+1})^*$ defined by the above conormals ${\bf n}_i$
and the following $k_i$
$$k_1=\dots=k_p=k_{p+2}=0,\; k_{p+1}=\tau,\; k_{p+3}=\lambda.$$

The polytope $\Delta$ is a $\Delta_p$ bundle on $\Delta_1$ (see
\cite{M-T3}). When $p=2$,   $\Delta$ is the prism whose base is
the triangle of vertices $(0,0,0),$ $(\tau,0,0)$ and $(0,\tau,0)$
and whose ceiling is the triangle determined by $(0,0,\lambda),$
$(\tau,0,\lambda+a_1\tau)$ and $(0,\tau,\lambda+a_2\tau)$

We assume that the above polytope $\Delta$ is a Delzant polytope.
The manifold (\ref{DefiM}) is in this case
$$M_{\Delta}=
 \{ z\in{\mathbb C}^{p+3}\,:\, \sum_{i=1}^{p+1}|z_i|^2=\tau/\pi,\;
\;-\sum_{j=1}^p a_j|z_j|^2+|z_{p+2}|^2+|z_{p+3}|^2=\lambda/\pi
    \}/ \simeq,$$ where $(z_j)\simeq (z'_j)$ iff there are
$\alpha,\beta\in U(1)$ such that
$$z'_j=\alpha\beta^{-a_j}z_j,\;j=1,\dots,p\,; \; \;\; z'_{p+1}=\alpha z_{p+1};\;\;\;z'_k=\beta z_k,\; k=p+2,\,p+3.$$

The symplectic form (\ref{omeg}) is
$$\omega_{\Delta}=(1/2)\big(\sigma_1+\dots+\sigma_p+\sigma_{p+2}),$$
where $\sigma_k= d\rho_k^2\wedge d\varphi_k.$

 And the moment map
\begin{equation}\label{momentb}
\mu_{\Delta}([z])=\pi(\rho_1^2,\dots, \rho_{p}^2,\,\rho_{p+2}^2).
\end{equation}

Thus $M_{\Delta}$  is the total space of the fibre bundle
${\mathbb P}(L_1\oplus\dots \oplus L_p\oplus{\mathbb C})\to
{\mathbb C}P^1$, where $L_j$ is the holomorphic line bundle over
${\mathbb C}P^1$ with Chern number $a_j$.

\smallskip

Given ${\bf \Hat b}=(b_1,\dots,b_p,0)\in{\mathbb Z}^{p+1}$ we
write
$$B:=\sum_{j=1}^pb_j,\;\;\;\;\; {\bf a}\cdot {\bf \Hat
b}:=\sum_{j=1}^pa_jb_j.$$

\begin{Prop}\label{CmDeltabundle}
 Let ${\bf \Hat b}=(b_1,\dots,b_p,0)$ be an element of
${\mathbb Z}^{p+1}$, then
 \begin{equation}\label{CmbMT}
\langle \rm{Cm}(\Delta),\,{\bf\Hat b}\rangle=\frac{\tau}{p+2}
\frac{\lambda(p+2)B+\tau (AB+{\bf a}\cdot{\bf \Hat
b})}{\lambda(p+1)+\tau A}.
\end{equation}

\end{Prop}

{\it Proof.} By Lemma \ref{Lemmaint}
 \begin{equation}\label{intMomega}
  \int_{M_{\Delta}}(\omega_{\Delta})^{p+1}=(p+1)!\int_{S_p(\tau)}\big(\lambda+\sum_{j=1}^p
a_jx_j\big)=(p+1)!\Big(\frac{\lambda\tau^p}{p\,!}+\frac{\tau^{p+1}A}{(p+1)!}\Big).
 \end{equation}
Similarly, for $k=1,\dots, p$
$$\int_{M_{\Delta}}x_k\,(\omega_{\Delta})^{p+1}=(p+1)!\Big(\frac{\lambda\tau^{p+1}}{(p+1)!}+
 \frac{\tau^{p+2}}{(p+2)!}\sum_{j\ne k}a_j+  \frac{2\tau^{p+2}a_k}{(p+2)!}
 \Big).$$
   So  the $k$-th coordinate of $\text{Cm}(\Delta)$, $\bar{x_k}$,
  is
 $$\bar{x_k}=\frac{\tau}{p+2}\,  \frac{\lambda(p+2)+\tau\big(A+ a_k   \big)}
 {\lambda(p+1)+\tau A}.$$
 So
  $$\langle\text{Cm}(\Delta),\,{\bf\Hat
 b}\rangle=\frac{\tau}{p+2}\frac{\lambda(p+2)B+\tau (AB+{\bf a}\cdot{\bf \Hat
b})}{\lambda(p+1)+\tau A}.$$

\qed

\begin{Thm}\label{ThmHatb}
Let ${\bf \Hat b}=(b_1,\dots,b_p,0)$ be an element of ${\mathbb
Z}^{p+1}$. For all for all $\tau>0$  and  all $\lambda$ with
$\lambda+a_i>0$,
 $$I(\lambda,\tau;\,{\bf\Hat
b})={\mathcal K}(\lambda,\tau)\Hat{{\mathcal Z}},$$ where
\begin{equation}\label{MathcalK}{\mathcal
K}(\lambda,\tau)=\frac{\tau^{p+1}}{p+2}\Big(\frac{\tau}{\lambda(p+1)+\tau
A} -1\Big)\;\;\;\text{and}\;\;\;\Hat{{\mathcal Z}}=2((p+1){\bf
a}\cdot{\bf \Hat b}-AB).
 \end{equation}
\end{Thm}

 {\it Proof.}
 We write
 \begin{equation}\label{I(PhijPhi'j)}\frac{1}{p+1}I(\lambda,\tau;\,{\bf\Hat
 b})=\sum_{j=1}^{p+3}\big(\langle \text{Cm}(\Delta),\,{\bf\Hat
 b}\rangle\Phi_j -\Phi'_j   \big),
 \end{equation}
 where
$$\Phi_j:=\int_{z_j=0}(\omega_{\Delta})^p,\;\;\;\;\Phi'_j:=\int_{z_j=0}\langle\mu_{\Delta},{\bf
b}\rangle(\omega_{\Delta})^p.$$

To calculate the values of the $\Phi_i$, we will distinguish
  three    cases, according to the value of $i$:

  {\it a)} $i=1,\dots,p$

 {\it b)} $i=p+1$

 {\it c)} $i=p+2,\,p+3.$

We will also respect  this classification in the calculation of
the $\Phi'_i$.

\smallskip

  {\it a)}  We calculate the value of $\Phi_1$.  On $z_1=0$
$$(\omega_{\Delta})^p=\frac{p\,!}{2^p}\big(\sigma_2\wedge\dots\wedge{\sigma_{p}}\wedge\sigma_{p+2}\big).$$
If we put $x_i:=\pi\rho_i^2$, then
\begin{equation}\label{Phi1}
\frac{1}{p\,!}\Phi_1=\int_0^\tau dx_2 \int_0^{\tau-x_2} dx_3\dots
\int_0^{\tau -\sum_{j=2}^{p-1}x_j}dx_p \int_0^{\lambda
+\sum_{j=2}^{p}a_jx_j}dx_{p+2}.
\end{equation}
It follows from Lemma \ref{Lemmaint}
 $$\frac{1}{p\,!}\Phi_1=\int_{S_{p-1}(\tau)}\big(\lambda+\sum_{j=2}^{p}a_jx_j\big)=
 \frac{\lambda\tau^{p-1}}{(p-1)\,!}+\frac{\tau^p}{p\,!}\sum_{j=2}^pa_j.$$

For $k=1,\dots, p$, a similar calculation gives
 \begin{equation}\label{Phikhatb}
 \frac{1}{p\,!}\Phi_k=
 \frac{\lambda\tau^{p-1}}{(p-1)\,!}+\frac{\tau^p}{p\,!}\sum_{j\ne k}^pa_j.
  \end{equation}

{\it b)} Next we consider $\Phi_{p+1}$. Now $x_{p+1}=0$. So
$x_p=\tau-\sum_{j=1}^{p-1}x_j$ and
$$-\sum_{j=1}^{p-1}a_{jp}x_j-a_p\tau +x_{p+2}+x_{p+3}=\lambda,\;\; \text{with}\; a_{jp}:=a_j-a_p.$$
Hence
\begin{equation} \label{alignPhi}
\frac{1}{p\,!}\Phi_{p+1}
=\int_{S_{p-1}(\tau)}\big(\lambda+\sum_{j=1}^{p-1}a_{jp}x_j
+a_p\tau\big) =
 \frac{\lambda\tau^{p-1}}{(p-1)\,!}+  \frac{\tau^p}{p\,!}A.
  \end{equation}

{\it c)}

$ \Phi_{p+2}= \Phi_{p+3}= \tau^p$.

Thus, it follows from (\ref{alignPhi}) and (\ref{Phikhatb})
\begin{equation}\label{sunPhikHatb}
\sum_{j=1}^{p+3}\Phi_j=(2+pA)\tau^p+p(p+1)\lambda\tau^{p-1}
\end{equation}

\smallskip

Next we determine the values of the $\Phi'_i$.

 {\it a')} On $z_k=0$, with $k=1,\dots, p$, $\,\langle\mu_{\Delta}[z],\,{\bf\Hat
 b}\rangle=\sum_{i\ne k}^{p}b_ix_i$. Then
$$\frac{1}{p\,!}\Phi'_k=\int_{S_{p-1}(\tau)}\big(\lambda+\sum_{j\ne k}^{p}a_jx_j\big)\sum_{i\ne k}^{p}b_ix_i.$$
By Lemma \ref{Lemmaint}
 \begin{equation}\label{Phi'1Hat}
 \frac{1}{p\,!}\Phi'_k=\Big(\sum_{i\ne k}^{p}b_i\Big) \frac{\lambda\tau^
p}{p\,!}+ 2\Big(\sum_{i\ne k}a_ib_i \Big)
\frac{\tau^{p+1}}{(p+1)\,!}
 +\Big(\sum_{j\ne i,\,j\ne k\ne i}a_ib_j\Big)  \frac{\tau^{p+1}}{(p+1)\,!}
\end{equation}

{\it b')} Using the notation introduced in {\it b)},
 $$\frac{1}{p\,!}\Phi'_{p+1}=\int_{S_{p-1}(\tau)}\Big(\lambda+\sum_{j=1}^{p-1}a_{jp}x_j
+a_p\tau\Big)\Big(\sum_{i=1}^{p-1}b_{ip}x_i+b_p\tau\Big),$$ with
$b_{ip}:=b_i-b_p$.
 Lemma \ref{Lemmaint} and a straightforward calculation give
 \begin{equation}\label{Phi'Hatb}
  \Phi'_{p+1}=\lambda\tau^pB+\big({\bf a}\cdot {\bf\Hat b}+AB
 \big)\frac{\tau^{p+1}}{p+1}
 \end{equation}

{\it c')} $$\Phi'_{p+2}=\Phi'_{p+3}=\frac{B\tau^{p+1}}{p+1}.$$

From this last result together with (\ref{Phi'1Hat}) and
(\ref{Phi'Hatb}) it follows
\begin{equation}\label{SumPhi'Hat}
\sum_{j=1}^{p+3}\Phi'_j=B\lambda\tau^p+\Big((p+1){\bf
a}\cdot{\bf\Hat b}+(p-1)AB+2B\Big)\frac{\tau^{p+1}}{p+1}
\end{equation}

Taking into account (\ref{I(PhijPhi'j)}), Proposition
\ref{CmDeltabundle}, (\ref{SumPhi'Hat}) and (\ref{sunPhikHatb}),
by means of an easy but tedious calculation   one obtains
 \begin{equation}\label{mathcalHatZ}
 I(\Delta,\,{\bf\Hat b})=\frac{2((p+1){\bf a}\cdot {\bf
b}-AB)}{(p+2)(\lambda(p+1)+\tau
A)}\big((1-A)\tau^{p+2}-(p+1)\lambda\tau^{p+1}  \big).
\end{equation}
If we define
 $\Hat{\mathcal Z}$ and   ${\mathcal K}(\tau,\lambda)$ as in the
 statement of theorem, then (\ref{mathcalHatZ}) can be written
 $I(\Delta,\,{\bf\Hat b})={\mathcal K}(\tau,\lambda)\Hat{\mathcal
 Z}.$

 \qed

We write ${\bf \Dot b}$ for the element $(0,\dots,0,b)\in{\mathbb
Z}^{p+1}$.

\begin{Prop}\label{CmDeltabundleDot}
 Given ${\bf \Dot b}=(0,\dots,0,b)\in{\mathbb
Z}^{p+1}$, then
 \begin{equation}\label{CmbMTDot}
\langle \rm{Cm}(\Delta),\,{\bf\Dot b}\rangle=
\frac{b}{2}\,\frac{(p+1)(p+2)\lambda^2+2(p+2)A\lambda\tau+({\bf
a}\cdot {\bf a}+A^2)\tau^2}{(p+2)\big((p+1)\lambda+A \tau \big)},
\end{equation}
where ${\bf a}\cdot {\bf a}=\sum_ia_i^2$.
\end{Prop}

{\it Proof.} We need to calculate $\int_M bx_{p+2}\omega^{p+1}.$
By Lemma \ref{Lemmaint}
 \begin{align} \label{alignCM}
  \frac{1}{(p+1)!}\int_M
bx_{p+2}\omega^{p+1}=&\frac{b}{2}\int_{S_p(\tau)}\big(\lambda+\sum_{j=1}^pa_jx_j
\big)^2 \\ \notag
 =&\frac{b}{2}\Big(\frac{\lambda^2\tau^p}{p\,!}+
\frac{2A\lambda\tau^{p+1}}{(p+1)!}+
 \frac{({\bf a}\cdot{\bf a}+A^2)\tau^{p+2}}{(p+2)!} \Big).
 \end{align}
Formula (\ref{CmbMTDot}) is a consequence of (\ref{intMomega})
together with (\ref{alignCM}).

\qed

\begin{Thm}\label{ThmIDot}
Let ${\bf \Dot b}=(0,\dots,0,b)$ be an element of ${\mathbb
Z}^{p+1}$. For all for all $\tau>0$  and  all $\lambda$ with
$\lambda+a_i>0$,
 $$I(\lambda,\tau;\,{\bf\Dot
b})={\mathcal K}(\lambda,\tau)\Dot{{\mathcal Z}},$$ where
 $$\Dot{{\mathcal Z}}=b\big((p+1){\bf a}\cdot{\bf
 a}-A^2\big),$$
 and ${\mathcal K}(\lambda,\tau)$ is defined as in
 (\ref{MathcalK}).
 \end{Thm}

{\it Proof.} As in the preceding theorem
\begin{equation}\label{I(PhijPhi'jDot)}
\frac{1}{p+1}I(\lambda,\tau;\,{\bf\Dot
 b})=\sum_{j=1}^{p+3}\big(\langle \text{Cm}(\Delta),\,{\bf\Hat
 b}\rangle\Phi_j -\Phi'_j   \big).
 \end{equation}
  The
$\Phi_j$ in (\ref{I(PhijPhi'jDot)}) as the same as in
(\ref{I(PhijPhi'j)}). But now $\langle\mu,\,{\bf\Dot
b}\rangle=bx_{p+2}$. We will follow the same steps as in the proof
of Theorem \ref{ThmHatb} for calculating the $\Phi'_j$.

{\it a')}
 The corresponding
$\Phi'_k=\int_M\langle\mu,\,{\bf\Dot b}\rangle\omega^{p}$ can be
calculated as in Theorem \ref{ThmHatb}.
 For $k=1,\dots,p$ one has
 $$\frac{1}{p\,!}\Phi'_k=\frac{b}{2}\Big( \frac{\lambda^2\tau^
{p-1}}{(p-1)\,!}
 +2\frac{\lambda\tau^{p}}{p\,!}\sum_{j\ne k}a_j
 +2\frac{\tau^{p+1}}{(p+1)\,!}\sum_{j\ne k}a_j^2
  + \frac{\tau^{p+1}}{(p+1)\,!} \sum_{j\ne i,\,j\ne k\ne i}a_ia_j
 \Big). $$

{\it b')} Now
  $$\frac{1}{p\,!}\Phi'_{p+1}=\frac{b}{2}\Big(
  \frac{\lambda^2\tau^{p-1}}{(p-1)!}+
  \frac{2\lambda\tau^p A}{p\,!}+
   \frac{\big({\bf a}\cdot {\bf a}+A^2\big)\tau^{p+1}}{(p+1)!}
  \Big)$$

{\it c')} In this case
$$\Phi'_{p+2}=0,\;\;\;\; \frac{1}{p\,!}\Phi'_{p+3}=b\Big( \frac{\lambda\tau^p}{p\,!}   +
\frac{A\tau^{p+1}}{(p+1)!}\Big).$$

If we insert in (\ref{I(PhijPhi'jDot)}) the values for the
$\Phi_j$ obtained in the proof of Theorem \ref{ThmHatb}, the above
values of the $\Phi'_j$ and   (\ref{CmbMTDot}) we arrive to
$I(\tau,\lambda;\,{\bf b})={\mathcal
K}(\lambda,\tau){\Dot{\mathcal Z}}.$

\qed

Given ${\bf b}=(b_1,\dots,b_p,b)\in{\mathbb Z}^{p+1}$, we write
${\bf \Hat b}=(b_1,\dots,b_p,0)$ and ${\bf \Dot b}:={\bf b}-{\bf
\Hat b}$. We put
 \begin{equation}\label{mathcalZ}
 {\mathcal Z}({\bf b}):=(p+1)\big({\bf a}\cdot(2{\bf \Hat b}+b{\bf
a})\big)-A(2B+bA),
 \end{equation}
 that is,
 ${\mathcal Z}={\Hat{\mathcal Z}}+{\Dot{\mathcal Z}}.$
It follows from (\ref{Ib}) that $I(\tau,\lambda;\,{\bf b})$ is a
group homomorphism with respect to the variable ${\bf b}$. By
 Theorem \ref{ThmHatb} and Theorem  \ref{ThmIDot}  one has

\begin{Thm}\label{ThmHatb+Dotb}
If ${\bf b}=(b_1,\dots,b_{p},b)\in{\mathbb Z}^{p+1}$, then
$$I(\tau,\lambda;\,{\bf b})={\mathcal
K}(\lambda,\tau){\mathcal Z}({\bf b}),$$
 where ${\mathcal Z}({\bf b})$ is given by (\ref{mathcalZ}).
 \end{Thm}

 This theorem expresses  the value of $I(\tau,\lambda;\,{\bf b})$ as
 the product of two factors. ${\mathcal K}$ is independent of the
 Hamiltonian loop, it depends only of $\lambda,\tau$. On the contrary, the factor
 ${\mathcal Z}$ is constant on the chamber of the polytope.

\smallskip

 Let ${\bf b}={\bf \Hat b}+{\bf \Dot b}$ be  as  before, since
 $\langle\rm{Cm}(\Delta),\,{\bf
b}\rangle=\langle\rm{Cm}(\Delta),\,{\bf\Hat b}\rangle
+\langle\rm{Cm}(\Delta),\,{\bf\Dot b}\rangle$, by (\ref{CmbMT})
and (\ref{CmbMTDot})
$$\lim_{\lambda\to 0}\langle  \rm{Cm}(\Delta),\,{\bf b}\rangle=
\frac{b}{2}\,\frac{ ({\bf a}\cdot{\bf a}+A^2  )\tau} {(p+2)A}+
\frac{ ({\bf a}\cdot{\bf \Hat b}+AB  )\tau} {(p+2)A}.$$
$$\lim_{\tau\to 0}\langle  \rm{Cm}(\Delta),\,{\bf b}\rangle=
\frac{b\lambda}{2}.$$
 Therefore,  $(\Delta,\,{\bf b})$ is a mass
iff
 \begin{equation}\label{Masslinear}\langle\rm{Cm}(\Delta),\,{\bf b}\rangle=\frac{b\lambda}{2}+
\Big( \frac{b}{2}\,\frac{ ({\bf a}\cdot{\bf a}+A^2  )} {(p+2)A}+
\frac{ ({\bf a}\cdot{\bf \Hat b}+AB  )} {(p+2)A}\Big)\tau
 \end{equation}
 If we insert in the equation
(\ref{Masslinear}) the expressions for $\langle
\rm{Cm}(\Delta),\,{\bf \Hat b}\rangle$ and $\langle
\rm{Cm}(\Delta),\,{\bf \Dot b}\rangle$ given by (\ref{CmbMT}) and
(\ref{CmbMTDot}) we obtain the following equivalent condition for
the pair $(\Delta,\,{\bf b})$ to be mass linear.

\begin{equation}\label{equvMassLinear}
B(p+2)+bA(p+2)=\frac{bA(p+2)}{2}+\Big( \frac{b({\bf a}\cdot{\bf
a}+A^2 )}{2A}+ \frac{AB+{\bf a}\cdot{\bf \Hat b}}{A}  \Big)(p+1).
\end{equation}
 That is,
 $$\frac{b( A^2-(p+1){\bf a}\cdot {\bf a} )}{2}=(p+1){\bf a}\cdot
 {\bf\Hat b}-AB.
 $$
Taking into account the definition
 of ${\mathcal Z}$ given in (\ref{mathcalZ})
 we can state the following theorem
\begin{Thm}\label{ThmMasLinear}
If ${\bf b}\in{\mathbb Z}^{p+1}$, then $(\Delta,\,{\bf b})$ is a
mass linear pair iff  ${\mathcal Z}({\bf b})=0$.
\end{Thm}

A consequence of Theorem \ref{ThmHatb+Dotb} and Theorem
\ref{ThmMasLinear} is the following result

\begin{Thm}\label{vanish}
 Assume that the Delzant polytope  $\Delta$ is a  $\Delta_p$ bundle over $\Delta_1$.
  Given ${\bf b}\in{\mathbb Z}^{p+1}$, then
the pair $(\Delta,{\bf b})$  is mass linear iff the
 the characteristic number $I(k;\,{\bf b})$ of the Hamiltonian
fibration $E_{\Delta(k),{\bf b}}\to S^2$   is zero for all $k$ in
the chamber of $\Delta$.
\end{Thm}

\smallskip

Given $\lambda$ and $\tau$, the map
 $\,{\bf b}\in{\mathbb Z}^{p+1}\mapsto I(\lambda,\tau;\,{\bf b})\in{\mathbb
 R}\,$ is a group homomorphism. By Theorem \ref{ThmHatb+Dotb} its
 kernel is ${\mathcal H}:=\{ {\bf b}\,|\,{\mathcal Z}({\bf b})=0  \}.$

Taking into account (\ref{Ihomo}), if ${\bf b}\notin {\mathcal H}$
then
 $1\ne[\psi_{\bf
 b}]\in\pi_1(\rm{Ham}(M_{\Delta}))$.
 So we have the following
 Proposition
 \begin{Prop}\label{Proppi}
 If $(\Delta,\,{\bf b})$ is not mass linear, then $\psi_{\bf b}$
 generates an infinite cyclic subgroup in
 $\pi_1(\rm{Ham}(M_{\Delta}))$.
\end{Prop}

In particular ${\bf\Hat b}=(b_1,\dots,b_p,0)$ belongs to
${\mathcal H}$ iff
 \begin{equation}\label{sumbiai}
 \sum_{j=1}^p
 b_j\big((p+1)a_j-\sum_{i=1}^p a_i\big)=0.
 \end{equation}

We put $c_j:=(p+1)a_j-\sum_{i=1}^pa_i$. If $(a_1,\dots,a_p)\ne 0$,
then ${\bf\Hat c}:=(c_1,\dots,c_p,0)$ does not satisfy
(\ref{sumbiai}).
 Thus one has the
following corollary
 \begin{Cor}\label{Coropi} If $(a_1,\dots,a_p)\ne (0,\dots,0)$, then $[\psi_{\bf\Hat c}]$ generates an infinite
 cyclic subgroup of $\pi_1(\rm{Ham}(M_{\Delta}))$.
  \end{Cor}

Analogously
\begin{Cor}\label{CoroDotb} If $(a_1,\dots,a_p)\ne (0,\dots,0)$ and ${\bf \Dot b}=(0,\dots,0,b)\ne 0$, then
$[\psi_{\bf\Dot b}]$ generates an infinite
 cyclic subgroup of $\pi_1(\rm{Ham}(M_{\Delta}))$.
  \end{Cor}

\smallskip


\section{One point blow up of ${\mathbb C}P^n$.}

\subsection{Hirzebruch surfaces.} \label{SubsectHirz}
Given $k\in{\Bbb Z}_{>0}$ and $\tau,\lambda\in{\Bbb R}_{>0}$ with
$\sigma:=\tau-k\lambda>0$, we consider the Hirzebruch surface $M$
determined by these numbers. $M$ is the quotient
$$\{z\in {\Bbb C}^4\,\,:\,\, |z_1|^2+k|z_2|^2+|z_4|^2=\tau/\pi,\;
|z_2|^2+|z_3|^2=\lambda/\pi\}/{\Bbb T},$$ where the equivalence
defined by ${\Bbb T}=(S^1)^2$ is given by
$$(a,b)\cdot(z_1,z_2,z_3,z_4)=(az_1, a^kbz_2,bz_3,az_4),$$
for $(a,b)\in (S^1)^2$. (The definition of Hirzebruch surface
given in \cite{V} can be obtained exchanging $z_1$ for $z_2$ in
the above definition.)

The manifold $M$ equipped with the following $(U(1))^2$ action
$$(\epsilon_1,\,\epsilon_2)[z_j]=[\epsilon_1 z_1,\,\epsilon_2
z_2,\,z_3,\,z_4],$$ is a toric manifold. The corresponding moment
 polytope $\Delta$ is the trapezium in ${\mathbb R}^2$ with
vertices $P_1=(0,0),$ $P_2=(0,\lambda),$ $P_3=(\tau,0),$
$P_4=(\sigma,\lambda).$ The mass center of $\Delta$ is
\begin{equation}\label{Cmblowup}{\rm Cm}(\Delta)=\Big(\frac{3\tau^2-3k\tau\lambda+k^2\lambda^2}{3(2\tau-k\lambda)},\,
\frac{3\lambda\tau-2k\lambda^2}{3(2\tau-k\lambda)} \Big).
\end{equation}
Given ${\bf b}=(b_1,\,b_2)\in{\mathbb Z}^2$, the pair $(\Delta,
\,{\bf b})$ is mass linear, iff there exist $A,B,C\in{\mathbb R}$
such that
$$\langle {\rm Cm}(\Delta),\,{\bf b}\rangle=A\tau+B\lambda+C,$$
for all $(\tau,\lambda)$ in the chamber of $\Delta$. From this
condition one deduces the following proposition
\begin{Prop}\label{linearHirz}
A necessary and sufficient condition for $(\Delta,{\bf b})$ to be
a mass linear pair is $2b_2=kb_1.$
\end{Prop}

{\it Remark.} The condition $2b_2=kb_1$ can be expressed in terms
of concepts introduced in \cite{M-T3}. The facets $P_1P_2$ and
$P_3P_4$ of $\Delta$ are equivalent (according to Definition 1.11
of \cite{M-T3}). Thus, if $2b_2=kb_1$, then ${\bf b}$ is
inessential (see Definition 1.13 in \cite{M-T3}).

\smallskip

 We denote by $\phi_t$ the following isotopy of $M$
$$\phi_t[z]=[e^{2\pi it}z_1,z_2,z_3,z_4].$$
$\phi$ is a loop in the Hamiltonian group of $M$. By $\phi'$ we
denote the Hamiltonian loop
$$\phi'_t[z]=[z_1,e^{2\pi it}z_2,z_3,z_4].$$
In Theorem 8 of \cite{V} we proved that
$I({\phi'})=(-2/k)I(\phi)$. If  ${\bf b}=(b_1,\,b_2)\in{\mathbb
Z}^2$, then
$$I({\psi_{\bf b}})=b_1I({\phi})+b_2I({\phi'})=(b_1-(2/k)b_2)I({\phi}).$$
From Proposition \ref{linearHirz} one deduces the following
theorem
\begin{Thm}\label{ThmHirz}
The pair $(\Delta,{\bf b})$ is mass linear iff $I({\psi_{\bf
b}})=0$. (Equivalently  $I(\tau,\lambda;\,{\bf b})=0$ for all
$(\tau,\lambda)$ in the chamber of $\Delta$.)
\end{Thm}

\smallskip

{\it Remark.} We will  deduce the vanishing of
 $I(\tau,\lambda;\,{\bf b})$ on ${\mathcal C}_{\Delta}$ when
$(\Delta,{\bf b})$ is mass linear, by an indirect way; that is,
without calculating the integrals of (\ref{Ib}).

 If $(\Delta,{\bf
b})$ is a mass linear pair, $I(\tau,\lambda;\,{\bf b})$ is  a
homogeneous polynomial of degree $2$ in $\tau,\lambda$, by
Proposition \ref{polyno}. That is,
\begin{equation}\label{C_1C_2C_3}
 I(\tau,\lambda;{\bf b})=C_1\lambda^2+C_2\lambda\tau+C_3\tau^2.
\end{equation}
Fixed $\tau$, if $\lambda\to 0$, then $\Delta$ converts into the
segment $[0,\,\tau]$, and ${\rm lim}_{\lambda\to 0}\,{\rm
Cm}(\Delta)=(\tau/2,\,0)$. In the limit, the facets of $\Delta$
give rise to the segments
$$F_1=[0,\,\tau],\;F_2=[0,\,\sigma],\;F_3=[\sigma,\,\tau],$$
on the  axis of abscissas. So ${\rm lim}_{\lambda\to
0}\,I(\tau,\lambda;\,{\bf b})$ is the sum of the contributions of
  $D_j=\mu^{-1}(F_j),\,$ $j=1,2,3$ (see (\ref{Ibr})). As $F_2,F_3$ is
  a decomposition of $F_1$, then
  $$\int_{D_1}f\omega_{\Delta}=\sum_{j=2}^3\int_{D_j}f\omega_{\Delta},$$
  and ${\rm lim}_{\lambda\to
0}\,I(\tau,\lambda;\,{\bf b})=-4\int_{D_1}f\omega_{\Delta}.$ On
the other hand
$$\int_{D_1}f\omega_{\Delta}=\int_0^{\tau}xb_1dx-\frac{\tau
b_1}{2}\tau=0.$$
 Hence ${\rm lim}_{\lambda\to
0}\,I(\tau,\lambda;\,{\bf b})=0,$  and $C_3$ in (\ref{C_1C_2C_3})
is zero. This result can also be deduced from Proposition
\ref{PropDivisible}.

Now let us assume that $k=1$. We denote by ${\Delta'}(\tau)$ the
triangle of vertices $(0,0),$ $(0,\tau),$ and $(\tau,0)$; the
corresponding toric manifold is ${\mathbb C}P^2$. As ${\rm
Ham}({\mathbb C}P^2)$ has the homotopy type of $PU(3)$,  the group
homomorphism $I:\pi_1({\rm Ham}({\mathbb C}P^2))\to{\mathbb R}$ is
zero.

For $\lambda>>1$ and $0<\sigma<<1$, the difference between the
polytopes  $\Delta(\tau,\lambda)$ and ${\Delta'}(\tau)$ is a
triangle with small sides.
 Thus, by formula (\ref{Ib}) for the
characteristic class $I$, the expression $|I(\tau,\lambda;\,{\bf
b})-I(\Delta'(\tau);\,{\bf b})|$ will be as small as we wish, if
$\lambda$ is big enough and $\sigma$ is sufficiently small.
 As $I(\Delta'(\tau);\,{\bf b})=0$, we conclude that
 $${\rm lim}_{\lambda\to
 \infty;\,\sigma\to 0}\,I(\lambda,\tau;\,{\bf b})=0.$$
If we take $\sigma=a/\lambda$, with $a$ an arbitrary positive
number, then
$$0={\rm lim}_{\lambda\to \infty}\,\big((C_1+C_2)\lambda^2+aC_2
\big).$$
 That is, $C_i=0$; and  by (\ref{C_1C_2C_3}), $I(\tau,\lambda;{\bf b})=0$ on the chamber of $\Delta$.


\subsection{One point blow up of ${\mathbb C}P^n$.} \label{SubsectCP2}

 The mass center of the simplex
$S_n(\tau)$, defined in Lemma \ref{Lemmaint}, is the point
 \begin{equation}\label{CmSn}{\rm
Cm}(S_n(\tau))=\frac{\tau}{n+1}w,
 \end{equation}
  with $w=(1,\dots,1)$. So any
pair $(S_n(\tau),\,{\bf b})$ with ${\bf b}\in {\mathbb Z}^n$  is
mass linear. Moreover, if we write
$$(1/n) \,I(S_n(\tau);\,{\bf b})=\langle{\rm
Cm}(S_n(\tau)),\,{\bf
b}\rangle\sum_{j=1}^{n+1}\Phi_j-\sum_{j=1}^{n+1}\Phi'_j,$$
  here
$$\Phi_j= \tau^{n-1},\;\;\; j=1,\dots, n+1.$$
$$\Phi'_j=\frac{\tau^{n}}{n}\sum_{i\ne j}b_i,\;\; j=1,\dots,
n.\;\;\; \;\;\; \;  \Phi'_{n+1}=\frac{\tau^n}{n}\sum_{j=1}^nb_j.$$
 Hence $$\langle{\rm
Cm}(S_n(\tau)),\,{\bf b}\rangle \sum_{j=1}^{n+1}\Phi_j=
 \tau^n \sum_{i=1}^nb_i =\sum_{j=1}^{n+1}\Phi'_j.$$
 We have the following proposition
 \begin{Prop}\label{CPn}
The invariant $I(S_{n}(\tau);\,{\bf b})$ is zero for all $\tau$
and all ${\bf b}\in{\mathbb Z}$.
\end{Prop}

 In this subsection $\Delta$ will be
\begin{equation}\label{Deltaonepoint}
\Delta=\Big\{(x_1,\dots,x_n)\in{\mathbb R}^n\,|\,
\sum_{i=1}^nx_i\leq\tau,\;0\leq x_i,\; x_n\leq\lambda\Big\},
\end{equation}
 where $\tau,\lambda\in{\mathbb R}_{>0}$ and
$\sigma:=\tau-\lambda>0$. That is, $\Delta$ is the polytope
obtained truncating the simplex $S_n(\tau)$ by a ``horizontal"
hyperplane through the point $(0,\dots,0,\lambda).$

As the volume   of $S_n(\tau)$ is $\tau^n/n!$, it follows from
(\ref{CmSn})
$$(\tau^n-\sigma^n)\,{\rm
Cm}(\Delta)=\tau^n\frac{\tau}{n+1}w-\sigma^n\big(\frac{\sigma}{n+1}w+\lambda\,
e_n \big).$$ That is,
 \begin{equation}\label{Cmtruncated}
  {\rm Cm}(\Delta)=
\frac{1}{\tau^n-\sigma^n}\Big(\big(\frac{\tau^{n+1}-\sigma^{n+1}}{n+1}
\big)w-\lambda\sigma^ne_n \Big).
\end{equation}

The pair $(\Delta,\, {\bf b}=(b_1,\dots, b_n))$ is mass linear iff
there exist $A,B,C\in{\mathbb R}$ such that
$$\sum_{j=1}^{n-1}b_j\frac{\tau^{n+1}-\sigma^{n+1}}{n+1}+b_n\Big(
\frac{\tau^{n+1}-\sigma^{n+1}}{n+1}-(\tau-\sigma)\sigma^n
 \Big)=(A\tau+B\sigma+C)(\tau^n-\sigma^n),$$
 for all $\tau,\sigma$ ``admissible". A straightforward calculation proves  the
 following proposition
 \begin{Prop}\label{PropLinear}
The pair $(\Delta, {\bf b})$ is mass linear iff
$$b_n=\frac{1}{n}\sum_{j=1}^{n-1}b_j.$$
\end{Prop}

The manifold $M_{\Delta}$ associated with $\Delta$ is the one
point blow up of
 ${\mathbb C}P^n$. On the other hand,   the
 general arguments showed in
the  Remark after Theorem \ref{ThmHirz}
  allow us to  prove following
Proposition

 \begin{Prop}\label{Propfiniteorder}
Let  $\Delta$ be the polytope (\ref{Deltaonepoint}), if
$(\Delta,{\bf b})$ is a mass linear pair, then
$I(\tau,\lambda;\,{\bf b})=0$ on ${\mathcal C}_{\Delta}$.
\end{Prop}

{\it Proof.} By Proposition \ref{polyno}, $I(\tau,\lambda;\,{\bf
b})=\sum_{i=0}^nB_i\lambda^{n-i}\tau^i$.
From (\ref{Cmtruncated}) one obtains
$${\rm lim}_{\lambda\to 0}\,{\rm
Cm}(\Delta)=\frac{\tau}{n}\Hat w,$$ with $\Hat w=(1,\dots,1,0).$
Hence the contribution of the base $F:=\{x\in \Delta\,|\,x_n=0\}$
to ${\rm lim}_{\lambda\to 0}\,I(\tau,\lambda;\,{\bf b})$ is
proportional to
\begin{equation}\label{contributionD}
\Big(\frac{\tau}{n}\sum_{j=1}^{n-1}b_j\int_{S_{n-1}(\tau)}1-\int_{S_{n-1}(\tau)}\sum_{i=1}^{n-1}b_ix_i\Big)=0.
\end{equation}

On the other hand, the other facets $F_i$ of $\Delta$ different
from the base $F$ give rise to a decomposition of $F$  in the
limit $\lambda\to 0$. Thus, with the notation of (\ref{Ibr}),
$${\rm lim}_{\lambda\to 0}\sum_{D_i\ne D}\int_{D_i}f\omega^{n-1}=
{\rm lim}_{\lambda\to 0} \int_{D}f\omega^{n-1},$$ where
$D_i:=\mu^{-1}(F_i)$ and $D:=\mu^{-1}(F).$
 By (\ref{contributionD}) this limit vanishes. It follows from
 (\ref{Ib}) that ${\rm lim}_{\lambda\to 0}\,I(\tau,\lambda;\,{\bf
 b})=0$. So
 $$I(\tau,\lambda;\,{\bf
b})=\sum_{i=0}^{n-1}B_i\lambda^{n-i}\tau^i.$$

For $\lambda>>1$ and $0<\sigma<<1$, the difference between the
polytopes  $\Delta(\tau,\lambda)$ and $S_{n}(\tau)$
 is a polytope with small edges.
 Thus,
$|I(\tau,\lambda;\,{\bf b})-I(S_{n}(\tau);\,{\bf b})|$ will be as
small as we wish, when $\lambda$ is big enough and $\sigma$ is
sufficiently small. By Proposition \ref{CPn} $I(S_{n}(\tau);\,{\bf
b})=0$ for all $\tau$. Hence,
 ${\rm lim}_{\lambda\to
 \infty;\,\sigma\to 0}\,I(\tau,\lambda;\,{\bf b})=0.$

 We take $\sigma=\lambda^{-1/n}$; from
 $$0=\lim_{\lambda\to
 \infty}\sum_{i=0}^{n-1}B_i\lambda^{n-i}(\lambda+\lambda^{-1/n})^i,$$
 one obtains a homogeneous system of $n$ linearly independent
 equations
 $$B_0+\dots+B_{n-1}=0,\;\;
 B_1+2B_2+\dots+(n-1)B_{n-1}=0,\;\;...\;\; , B_{n-1}=0,$$
 for the $n$ constants. So $B_j=0$, and $I(\tau,\lambda;\,{\bf
 b})=0$ on ${\mathcal C}_{\Delta}$.
 \qed

\smallskip

Next we will prove the reciprocal proposition of the preceding
one. We denote by $F_j$   the facet of $\Delta$ defined by
$x_j=0$, for $j=1,\dots,n$; $\,F_{n+1}$ will be   the ``ceiling"
$x_n=\lambda$ and $F_{n+2}$   the facet $x_1+\dots+x_n=\tau$.  We
write as above
 \begin{equation}\label{Itaulambda}
 \frac{1}{n} I(\tau,\lambda;\,{\bf b})=\langle{\rm
Cm}(\Delta),\,{\bf
b}\rangle\sum_{j=1}^{n+2}\Phi_j-\sum_{j=1}^{n+2}\Phi'_j,
 \end{equation}
 with $\Phi_j=\int_{F_j}1,\;\; \Phi'_j=\int_{F_j}\sum_ib_ix_i$.
 By Lemma \ref{Lemmaint}
$$ \Phi_j= \tau^{n-1}-\sigma^{n-1};\; \text{for}\;
 j=1,\dots,n-1,n+2.\; \;\;  \Phi_n= \tau^{n-1}.
 \;\;\; \Phi_{n+1}= \sigma^{n-1}.$$
 Thus we have
 \begin{equation}\label{sumPHif}
 \sum_{j=1}^{n+2}\Phi_j=(n+1)\tau^{n-1}+(1-n)\sigma^{n-1}.
 \end{equation}

 If ${\bf b}=\Hat{\bf b}=(b_1,\dots,b_{n-1},0)$, by (\ref{Cmtruncated})
 $$\langle{\rm Cm}(\Delta),\,\Hat{\bf
 b}\rangle=\frac{1}{n+1}\,\frac{\tau^{n+1}-\sigma^{n+1}}{\tau^n-\sigma^n}\,\sum_{j=1}^{n-1}b_j.$$
Similarly,
\begin{equation}\label{sumPHi'f}
\sum_{j=1}^{n+2}\Phi'_j=\frac{1}{n}\big(n\tau^{n}+(2-n)\sigma^{n}
 \big)\sum_{j=1}^{n-1}b_j
 \end{equation}

If we put $\epsilon:=\sigma/\tau$, it follows from
(\ref{Itaulambda}), (\ref{sumPHif}) and (\ref{sumPHi'f})
 \begin{equation}\label{Itaulambdaepsilon}
 \frac{1}{n!}I(\tau,\lambda;\,\Hat{\bf b})=
 -\frac{\tau^n}{n+1}\frac{\epsilon^{n-1}}{(n-2)!}\,\sum_{j=1}^{n-1}b_j+
 O(\epsilon^n).
 \end{equation}

\smallskip

Next we determine $I(\tau,\lambda;\,\Dot{\bf b})$, when $\Dot {\bf
b}=(0,\dots,0,b_n)$. Now
$$\sum_{j=1}^{n+2}\Phi'_j=b_n \tau^n \big( 1+(1-n)\epsilon^{n-1}+(n-2)\epsilon^n
\big)$$ and
$$\langle{\rm Cm}(\Delta),\,\Dot{\bf
b}\rangle=b_n\frac{\tau}{n+1}\big(1-\frac{n\epsilon^n}{n+1}
\big)+O(\epsilon^{n+1}).$$
  Hence
\begin{equation}\label{ItaulambdaepsilonDot}
 \frac{1}{n!}I(\tau,\lambda;\,\Dot{\bf b})=
 b_n\frac{n\tau^n}{n+1}\frac{\epsilon^{n-1}}{(n-2)!}+
 O(\epsilon^n).
 \end{equation}

Now, if ${\bf b}=(b_1,\dots,b_n)$, it follows from
(\ref{Itaulambdaepsilon}), (\ref{ItaulambdaepsilonDot})
\begin{equation}\label{ItaulambdaepsilonGeneral}
\frac{1}{n!}I(\tau,\lambda;\,{\bf b})=
 \Big(n\,b_n-\sum_{j=1}^{n-1}b_j\Big)   \frac{\tau^n}{n+1}\frac{\epsilon^{n-1}}{(n-2)!}+
 O(\epsilon^n).
\end{equation}
 From (\ref{ItaulambdaepsilonGeneral}) together with Proposition \ref{PropLinear} we deduce the
following proposition
\begin{Prop}\label{Propnecessbn=}
If $\Delta$ is the polytope (\ref{Deltaonepoint}) and $I(k;\,{\bf
b})=0$ for all $k\in{\mathcal C}_{\Delta}$, then $(\Delta,\,{\bf
b})$ is a mass linear pair.
\end{Prop}

Propositions \ref{Propfiniteorder}  and \ref{Propnecessbn=} imply
the following theorem
\begin{Thm}\label{Thmequiv}
If $\Delta$ is the polytope (\ref{Deltaonepoint}), then
 $(\Delta,\,{\bf b})$ is a mass linear pair iff
 $I(k;\,{\bf b})=0$ for all $k\in{\mathcal
 C}_{\Delta}$.
\end{Thm}

It follows from Theorem \ref{Thmequiv}, together with Proposition
\ref{PropLinear} and the homomorphism (\ref{Ihomo}) the following
proposition
 \begin{Prop}\label{Propblowup}
 If ${\bf b}=(b_1,\dots,b_n)\in{\mathbb Z}^n$ and
$\sum_{j=1}^{n-1}b_j\ne nb_n$, then $\psi_{\bf b}$ generates an
infinite cyclic subgroup in  $\pi_1(\rm{Ham}(M_{\Delta}))$.
\end{Prop}

\smallskip

{\it Remark.} When $n=3$ the   toric manifold $M$ corresponding to
$\Delta$ is
 $$M=\{z\in{\Bbb
C}^5\,:\,|z_1|^2+|z_2|^2+|z_3|^2+|z_5|^2=\tau/\pi,\,\,|z_3|^2+|z_4|^2=\lambda/\pi
\}/{\Bbb T},$$
 where the action of ${\Bbb T}=(U(1))^2$ is
defined by
\begin{equation}\label{ActD}
 (a,b)(z_1,z_2,z_3,z_4,z_5)=(az_1,az_2,abz_3,bz_4,az_5),
 \end{equation}
  for
$a,b\in U(1)$.

We consider the following loops in the Hamiltonian group of $M$
$$\psi_t[z]=[z_1e^{2\pi it},z_2,z_3,z_4,z_5],\;\;
\psi'_t[z]=[z_1,z_2e^{2\pi it},z_3,z_4,z_5],$$
$$\tilde\psi_t[z]=[z_1,z_2,z_3e^{2\pi it},z_4,z_5].$$
In \cite{V2} (Remark in Section 4) we gave formulas that relate
 the characteristic numbers
associated with these loops
$$I({\psi})=I({\psi'})=(-1/3)I({\tilde\psi}).$$
 So for ${\bf b}=(b_1,b_2,b_3)\in{\mathbb Z}^3$,
 \begin{equation}\label{RemarkIpsi}
 I({\psi_{\bf b}})=(b_1+b_2-3b_3)I({\psi}).
\end{equation}
 By Proposition
\ref{PropLinear} the vanishing of $I({\psi_{\bf b}})$ in
(\ref{RemarkIpsi}) is equivalent  to the fact that $(\Delta,\,{\bf
b})$ is a mass linear pair. This equivalence is a particular case
of Theorem \ref{Thmequiv}.


\end{document}